\documentclass[12pt]{article}

\usepackage{amsmath} 
\usepackage{graphicx} 
\usepackage{makeidx}	
\usepackage{latexsym}
\usepackage{amssymb}
\usepackage{amsfonts}
\newfont{\rams}{msbm10 scaled\magstep1} 
\newtheorem{teo}{Theorem}[section]
\newtheorem{cor}[teo]{Corollary}

\newtheorem{pro}[teo]{Proposition}

\newtheorem{defn}[teo]{Definition}
\newtheorem{rem}[teo]{Remark}

\newenvironment{dimo}%
	{\underline{\emph{Proof}} \newline }%
	{\hfill $\Box$ \newline }

\newenvironment{ex}%
	{\underline{\emph{Example}} \newline }%
	{\hfill $\Box$ \newline \newline}
	
\makeindex 

\begin{document}
\title{Equivalence of three different definitions of irreducible element} 
\author{Ornella Greco}
\date{\empty}
\maketitle
\begin{abstract} 
\footnotesize{
In this work, we try to draw a comparison among the three different concepts of irreducibility, given by Fletcher, Galovich, Bouvier; we have found out that these three concepts are equivalent in a particular class of rings, called 'rings with only harmless zero divisors'.}
\end{abstract}

\section{Rings with only harmless zerodivisors}

Galovich, Bouvier and Fletcher generalized the concept of unique factorization domain to rings with zero divisors: to do this, they gave new definitions of irreducibile elements and of associate elements (see \cite{and2}, \cite{and1} for an argumentation  about unique factorization rings).\\ 
The aim of this paper is to compare these three different definitions of irreducible element, and to prove their equivalence under certain conditions. Namely, in the class of ring that  we will now present the equivalence holds.\\

Let $R$ be a commutative ring with unity, let us denote by $Z(R)$ the set of all zero divisors in $R$, by $U(R)$ the group of units of $R$, by $Nil(R)$ the nilradical of $R$, by $J(R)$ its Jacobson radical.
\begin{defn}
		 Let $R$ be a commutative ring, we say that $r \in R$ is a \emph{harmless zero divisor} if $r\in Z(R)$ and there exists a unit $u$ such that $r=1-u$.\\
		  A ring $R$ is said to be a ring with only harmless zero divisors if every zero divisor in $R$ is harmless.
\end{defn} 
We notice that integral domains and local rings belong to this class of
rings: integral domains have not zero divisors so the property is trivially
fulfilled; also a local ring $(R,M)$ is a ring with only harmless zero divisors, since $Z(R) \subseteq M=J(M)\subseteq 1-U(R)$.\\
But there are rings with only harmless zero divisors that are neither local rings nor integral domains, as the one described in the following example.\\

\noindent
\begin{ex}
	Let $A$ be an Artinian, local, principal ideal ring, then $A[x]$ is  a ring with only harmless zero divisors.\\
If $(t)$ is the only maximal ideal of $A$, one can prove  that 
\[
	(t)=J(A[x])=Z(A[x]),
\] 
moreover in every commutative ring with unity, $R$, we have that $J(R)\subseteq 1-U(R)$, so  we have that $A[x]$ is a ring with only harmless zero divisors.\\
$A[x]$ is not an integral domain, since $Z(A[x])$ is not empty. Furthermore, $A[x]$ is not a local ring.
\end{ex}

\subsection{Different concepts of irreducibility}
In the following we present the three different definitions of irreducible element: the classical definition that Galovich (\cite{galovich}) adopted in his definition of unique factorization ring; the definition given by Bouvier in \cite{Bouvier}; the definition given by Fletcher in \cite{fletcher1}.\\
These definitions are equivalent in rings with only harmless zero divisors, but not in general.\\

\noindent
Let $R$ be a commutative ring with unity.
\begin{defn}
We say that $r \in R$ is  \emph{irreducible} if 
\[
r= ab \  \Rightarrow \ a\ \textrm{is a unit or }\ b\ \textrm{is a unit}.
\]
\end{defn}

\begin{defn}
			Let $r$ be a non-zero and non-unit element in $R$, we say that $r$ is \emph{B-irreducible} if the ideal $(r)$ is
			a maximal element in the set of all the principal proper ideals of $R$, ordered by the inclusion relation.

	\end{defn}
At first, we want to draw a comparison among these two definitions.
\begin{pro}\label{gb-irr}
			Let $r\in R$ be a non-unit, non-zero element, and suppose that $r$ is irreducible,
			then $r$ is B-irreducible.\\
			The converse is not true in general.
	
		\end{pro}
		\begin{dimo}
		By contradiction, let $(s)$ be a proper, principal ideal of $R$, and let $(r)\subsetneq (s)$.
		So there is $a\in R$, such that $r=as$, but $s $ is not a unit and $r$ is irreducible, hence
		$a$ must be a unit. It follows that $(r)=(s)$, against assumption.
		
		For the second part of the proof, we give an example of a ring, in which there is a B-irreducible element
		that is not irreducible. Let us consider $R= \mathbb{Z}_6$, and $r=\overline{3}$: the ideal $(\overline{3})$
		is maximal among the principal proper ideals of $R$, but $\overline{3}=\overline{3} \cdot \overline{3}$, so we have that
		$\overline{3}$ is B-irreducible, but it is not irreducible.
		\end{dimo}
		
		\noindent
	In the above proposition, we have found out that the concept of B-irreducible element is stronger than the one of irreducible element, but in rings with only harmless zero divisors, these two concepts are the same. 
		
		\begin{pro}\label{bg-irr}
		  Let $R$ be a ring with only harmless zero divisors, then every B-irreducible element is irreducible.
		\end{pro}
			\begin{dimo}
			 	By contradiction, suppose that there is a non-zero, non-unit element $x$ in $R$ that is B-irreducible, but not irreducible. Then, the principal ideal generated by $x$ is a maximal element in the set of the principal and proper ideals of $R$; on the other hand, there are two non-unit elements, $a,b \in R$, such that $x=ab$  Since $x$ is B-irreducible, we have that $(x)=(a)=(b)$, so we get the relation $x(xcd-1)=0$, for some $c,d$.  We now distinguish between two cases: first, if $x$ is not a zero divisor, then $xcd=1$, and so $x$ is a unit, here we have a contradiction; second, if $x$ is  a zero divisor, then $xcd$ is still a zero divisor, and, by hypothesis, $1-xcd$ is a unit, and $x=0$, a contradiction. 
			\end{dimo}


Let us present the definition given by Fletcher: to do this, we need first the concept of  a refinement of a factorization.
\begin{defn}
	Let $r=a_1\cdots a_n$ be a factorization of $r\in R$. A \emph{refinement} of this factorization is obtained by factoring one or more of the factors.
\end{defn}

\begin{defn}
	A non-unit element $r\in R$ is said to be \emph{F-irreducible} if each factorization of $r$ has a refinement containing $r$, as one of the new factors.
\end{defn}

\noindent
This definition can be formulated in another, more intuitive, way, because of the following simple result.
\begin{rem}
The following conditions are equivalent:
\begin{enumerate}
\item $r\in R$ is an F-irreducible element;
\item if $r=ab$, then $a\in (r)$ or $b\in (r)$;
\item if $r=ab$, then $(r)=(a)$ or $(r)=(b)$.
\end{enumerate}

\end{rem}

We want now to compare this new definition of F-irreducible element with the one given by Bouvier.
 
\begin{pro}\label{bf-irr}
 If $r\in R$ is a B-irreducible element, then $r$ is an F-irreducible element.
\end{pro}
	\begin{dimo}
	 Suppose that $r=a_1a_2 \cdots a_m$ and, for instance, let $a_1,a_2, \dots , a_s$ be non-units, hence $(r)\subseteq (a_i)$, for each $i=1, 2, \dots , s$, and, by hypothesis, we must have that $(r)=(a_i)$ for those $i$, i.e. there is a refinement of the given factorization that contains $r$, then, because we have taken an arbitrary factorization, $r$ is F-irreducible.
	\end{dimo}
	
	\noindent
	So, we have  proved that Fletcher's definition of irreducible element is a more general concept than Bouvier's one.
But, this two concepts are equivalent in the case of a ring with only harmless zero divisors. 
\begin{pro}\label{fb-irr}
Let $R$ be a ring with only harmless zero divisors, if $r\in R$ is a non-zero, F-irreducible element, then it is also a B-irreducible element.
\end{pro}
	\begin{dimo}
		By contradiction, let us suppose that $r\in R$ is F-irreducible, but B-reducible, then there exists a non unit $a\in R$ such that $(r)\subsetneq (a)$, i.e. there is  a non unit element $b \in R$ such that $r=ab$. Since $r$ is F-irreducible, we have that $a\in (r)$ or $b\in (r)$: in the first case, we get a contradiction, because we obtain that $(r)=(a)$; in the second case, we have that $b=cr$ and that $r(1-ac)=0$, but $r\neq 0 $, then we get that $1-ac=1-u$, where $u$ is a unit, and that $a$ is a unit, and here we have the contradiction.
	\end{dimo}

Using Proposition \ref{gb-irr}, Proposition \ref{bf-irr}, Proposition \ref{bg-irr} and Proposition \ref{gb-irr}, we obtain the following important result.
\begin{cor}
In a ring with only harmless zero divisors $R$, the three given concepts of irreducibility are equivalent.
\end{cor}
In the following, we will prove that also in the direct product of finitely many rings with only harmless zero divisors the three given concepts of irreducible element are equivalent.\\

To do this, we need to list some useful properties about the behavior of F-irreducible elements in the direct product of finitely many rings.

\begin{pro}\label{mfirred}
Let us consider the direct product of $n$ commutative rings with unity, $B=A_1 \oplus \cdots \oplus A_n$, and let $(a_1, \dots , a_n)\in B$ be a non-unit, non-zero element, if $(a_1, \dots , a_n)$ is an F-irreducible element, then $\exists \ i\in \{1, \dots , n\}$ such that $a_i$ is F-irreducible and $a_j$ is a unit for each $j\neq i$.
\end{pro}
	\begin{dimo}
			Let us consider the factorization
			\[
						(a_1,a_2,  \dots , a_n)=(a_1, 1,\dots , 1 )\cdots (1,1,\cdots , a_n),
			\]
			but $(a_1, \dots , a_n)$ is F-irreducible, so, by definition, this factorization has a refinement that contains $(a_1, \dots , a_n)$ as one of the new factors, i.e. there is $i\in \{1, \dots , n\}$, such that $(1, \dots, 1 , a_i , 1, \dots ,1)\in ((a_1, \dots , a_n))$. This means that $a_j$ is a unit for each $j\neq i$, and that $a_i$ is an F-irreducible element, since, if it had a factorization without refinement that contains it, we could easily find such a factorization for $(a_1, \dots , a_n)$, against the assumption of the F-irreducibility.
	\end{dimo}

\begin{pro}\label{equivalence}
Let us consider the direct product of finitely many, say $n$, rings with only harmless zero divisors, $B=A_1 \oplus \cdots \oplus A_n$, then in $B$ the three definition of irreducible element are equivalent.

\end{pro}
	\begin{dimo}
	We want to prove that if an element in $B$ is F-irreducible then it is irreducible; this would conclude the proof, since we already know that if an element is irreducible, it is B-irreducible, and that a B-irreducible element is F-irreducible.\\
	Let us consider an F-irreducible element, $r=(a_1, \dots , a_n)\in B$.	By Proposition \ref{mfirred}, $\exists \ i$ such that $a_i$ is F-irreducible and $a_j$ is a unit for each $j\neq i$. So, if $r=(\alpha_1, \dots , \alpha_n)(\beta_1, \dots , \beta_n)$, we must have that $\alpha_j, \beta_j$ are units for each $j\neq i$. Moreover, since $A_j$ is a ring with only harmless zero divisors, the element $a_i$ is also irreducible, so from $a_i=\alpha_i\beta_i$, we deduce that either $\alpha_i$ or $\beta_i$ is a unit: in the first case, $(\alpha_1, \dots , \alpha_n)$ is a unit in $B$, in the second one, $(\beta_1, \dots , \beta_n)$ is a unit. We have proved that $r$ is irreducible.
	
	\end{dimo}

Using a similar argument, it is possible to prove the following more general result.

\begin{pro}
Let $B$ be the direct product of a finite number, say $n$, of commutative rings with unity in which the three different definitions of irreducible element are equivalent, then the equivalence holds also in $B$.
\end{pro}

\noindent
{\bf Acknowledgments} I would like to thank Professor Ralf Fr\"oberg and Professor Christian Gottlieb for their important support.

\noindent
 {\scshape Royal Institute of Technology, Department of Mathematics, 10044 Stockholm, Sweden}.\\
 {\itshape E-mail address}: \texttt{ogreco@kth.se}

\clearpage
\addcontentsline{toc}{section}{Bibliography}
\begin{thebibliography}{9}
	\bibitem{and2} A. G. A\={g}arg\"{u}n, D. D. Anderson, S. Valdes-Leon,
		\emph{Unique factorization rings with zero divisors},
		Comm. Algebra, 27(4), pp. 1967-1974, 1999.
	
	\bibitem{and1} D.D. Anderson and R. Markanda,
		\emph{Unique factorization rings with zero divisors},
		Houston J. Math 11, corrigendum, pp. 15-30, 1985.
	\bibitem{Atiyah} M. F. Atiyah, I. G. MacDonald,
		\emph{Introduction to Commutative Algebra},
		Addison-Wesley Publishing Company, 1961.
	\bibitem{Bouvier}
 A. Bouvier,  
		\emph{Structure des anneaux \`{a} factorisation unique},	Publ. D\'{e}p. Math. (Lyon) 11, pp. 39-49, 1974.	
	\bibitem{fletcher1} C.R. Fletcher,
		\emph{Unique Factorization Rings},
		Proc. Cambridge Philos. Soc. 65, pp. 579-583, 1969.
	\bibitem{fletcher2} C.R. Fletcher,
		\emph{ The structure of unique factorization rings},
		Proc. Cambridge Philos. Soc. 67, pp. 535-540, 1970.	
	\bibitem{galovich} S. Galovich,
		\emph{Unique factorization rings with zero divisors},
		Mathematical Magazines 5, pp. 276-283, 1978.
	\bibitem{samuel}  P. Samuel, 
		\emph{Unique factorization}, Amer. Math. Monthly 75, pp. 945-952, 1968.		
	\bibitem{zar}   O. Zariski, P. Samuel,
		\emph{Commutative Algebra}, D. Van Norstrand Company, 1967.
	
	\end{thebibliography}
\end{document}